\documentclass[11pt]{amsart}
\usepackage{amsmath,amssymb,amscd,verbatim}
\usepackage{amssymb}
\usepackage{amsthm}
\usepackage{amsmath} \usepackage{latexsym}
\usepackage{color}%

\newcommand{\mc}{\mathcal}
\newcommand{\sub}{\subseteq}
\newcommand{\nsub}{\nsubseteq}

\newcommand{\lra}{\Leftrightarrow}

\newcommand{\ra}{\Rightarrow}

\newcommand{\sm}{\setminus}
\newcommand{\al}{\alpha}
\newcommand{\La}{\Lambda}
\newcommand{\be}{\beta}

\DeclareMathOperator{\rad}{rad}
\newcommand{\op}{\operatorname}
\newcommand{\tmax}{t\op{-Max}}

\newcommand{\Na}{R\langle X \rangle}
\newcommand{\tspec}{t\op{-Spec}}

\newtheorem{theorem}{Theorem}[section]
\newtheorem{lemma}[theorem]{Lemma}
\newtheorem{proposition}[theorem]{Proposition}
\newtheorem{corollary}[theorem]{Corollary}
\newtheorem{remark}[theorem]{Remark}
\theoremstyle{definition}

\begin{document}
    \title{Ring-theoretic properties of P$v$MDs}

\date{\today}



\author{Said El Baghdadi}

\address{Department of Mathematics, Facult\'{e} des Sciences et
Techniques, P.O. Box 523, Beni Mellal,\break Morocco}

\email{baghdadi@fstbm.ac.ma}

\author{Stefania Gabelli}

    \address{Dipartimento di Matematica, Universit\`{a} degli Studi  \textquoteleft\textquoteleft Roma Tre", Largo S.  L.  Murialdo, 1, 00146 Roma, Italy}

\email{gabelli@mat.uniroma3.it}

\begin{abstract}  We extend to Pr\"ufer $v$-multiplication domains some distinguished ring-theoretic properties of  Pr\"ufer domains. In particular we consider the $t\#\#$-property, the $t$-radical trace property, $w$-divisoriality and $w$-stability.
\end{abstract}

\maketitle

\section*{Introduction}

A Pr\"ufer $v$-multiplication domain, for short a P$v$MD, is a domain whose localizations at $t$-maximal ideals are valuation domains \cite{gr}. For this reason, the ideal-theoretic properties of valuation domains globalize to $t$-ideals of P$v$MDs and several properties of ideals of Pr\"ufer domains hold for $t$-ideals of P$v$MDs: for example a domain is a P$v$MD if and only if each $t$-finite $t$-ideal is $t$-invertible.
 The aim of this paper is to show that, introducing suitable $t$-analogues of some distinguished properties of integral domains,  Pr\"ufer domains and P$v$MDs have a similar behaviour also from the ring-theoretic point of view. We recall that the class of P$v$MDs, besides Pr\"ufer domains, includes Krull domains and GCD-domains.

 As a matter of fact, the $t$-operation is not always as good as the $w$-operation for extending certain properties that
hold in the classical case, that is in the $d$-operation setting.
Thus in general it is often more convenient to consider the
$w$-analogue of a given property (see for instance \cite{FM1, FM2, ss}). However in a P$v$MD the
$w$-operation and the $t$-operation coincide \cite{K} and one can
use indifferently these two star operations.

In Section 1, we deal with the $t\#\#$-property and the $tRTP$-property. P$v$MDs satisfying the $t\#\#$-property have been studied in \cite[Section 2]{GHL}.  Here we characterize P$v$MDs with the $tRTP$-property; getting for example that a P$v$MD is a $tRTP$-domain if and only if each $v$-finite divisorial ideal has at most finitely many minimal primes. Then, generalizing the Pr\"ufer case, we show that the $t\#\#$-property and the $tRTP$-property are strictly connected for P$v$MDs. Among other results, we prove that a P$v$MD satisfying the $t\#\#$-property is a $tRTP$-domain and that the converse holds if each $t$-prime is branched. We also show that an almost Krull domain satisfying the $t\#\#$-property is a Krull domain.

 In Section 2 we introduce the notion of $w$-stability and relate it to $w$-divisoriality, a property
 defined and studied in \cite{ss}.  First we show that the study of $w$-stability can be reduced to the $t$-local case. Then we use this result to extend to P$v$MDs some properties of stable and divisorial Pr\"ufer domains. For example, we prove that $w$-stability of $t$-primes enforces a P$v$MD to be a generalized Krull domains and that an integrally closed $w$-stable domain is precisely a  generalized Krull domain with $t$-finite character.  We also characterize $w$-stable $w$-divisorial P$v$MDs  and show that these domains behave like totally divisorial Pr\"ufer domains.

 \medskip
We assume that the reader is familiar with the language of star operations \cite[Sections 32 and 34]{g1}. We recall some definitions and basic properties.

Throughout this paper $R$ will denote an integral domain
with quotient field $K$ and we will assume that $R \neq K$.

A \emph{star operation} is
a map $I\to I^*$  from the set $F(R)$ of nonzero fractional ideals of $R$ to
itself such that:

(1) $R^* = R$ and $(aI)^* = aI^*$, for all $a \in K \smallsetminus
\{0\}$;

(2) $I \sub I^*$ and $I \sub J \ra I^* \sub J^*$;

(3) $I^{**} = I^*$.

 A star operation $*$ is of {\it finite type} if
 $I^* = \bigcup \{J^*\,;\, J \sub I$  and $J$ is fini\-te\-ly
ge\-ne\-ra\-ted\},
 for each $I \in F(R)$.
To any  star operation $*$, we can associate a star operation
$*_{f}$ of finite type
by defining $I^{*_{f}}= \bigcup J^*$, with the
union taken over all finitely generated ideals $J$ contained
in $I$. Clearly $I^{*_{f}} \sub I^*$.

 An ideal $I\in F(R)$ is  a \emph{$*$-ideal}  if $I= I^*$ and is $*$-{\it finite} if $ I^* =J^*$ for
 some finitely generated ideal $J$. A $*$-finite $*$-ideal is also called a \emph{$*$-ideal of finite type}.

 A $*$-{\it prime} is a prime $*$-ideal  and a $*$-{\it maximal} ideal is an ideal
that is maximal in the set of the proper
$*$-ideals. A $*$-maximal ideal  (if it exists) is a prime ideal.
 If $*$ is a star operation of finite type, an easy application of
Zorn's Lemma shows that the set $*$-Max$(R)$ of the $*$-maximal ideals of
$R$ is not empty. In this case,
 for each $I\in F(R)$, $I^* = \bigcap_{M\in
*\op{-Max}(R)}I^*R_{M}$; in particular
$R = \bigcap_{M\in *\op{-Max}(R)}R_{M}$ \cite{gr}.

For any star operation $*$, the
set of $*$-ideals of $R$ is a semigroup under the \emph{$*$-multiplication}, defined by
 $(I, J)\mapsto (IJ)^*$, with unity $R$. An ideal $I\in F(R)$ is called
\emph{$*$-invertible} if $I^*$ is invertible with respect to the
$*$-multiplication. In this case the $*$-{\it inverse} of
$I$ is $(R:I)$.
 Thus $I$ is  $*$-invertible if and only if $(I(R:I))^*=R$.

The identity is a star operation, called the \emph{$d$-operation}.
The \emph{$v$-operation} (or \emph{divisorial closure}), the
\emph{$t$-operation} and the \emph{$w$-operation} are the best known nontrivial star operations and are defined in the following way.
For each $I \in F(R)$, we set $I_v:=(R\colon (R\colon I))$ and $I_t:=\bigcup J_v$ with the
union taken over all finitely generated ideals $J$ contained
in $I$. Hence the $t$-operation is the finite type
star operation associated to the $v$-operation.
 The $w$-operation is the star operation of finite type defined by  setting $I_{w}:=
\bigcap_{M\in t\text{-Max}(R)}IR_{M}$.
We have $w$-Max$(R) = \tmax(R)$ and $IR_{M}= I_{w}R_{M}€\sub
I_{t}€R_{M}€$, for
each $M \in \tmax(R)$. Thus $I_w \sub I_t \sub I_v$.
Also, an ideal $I\in F(R)$ is $w$-invertible if and only if it is $t$-invertible.

A Pr\"ufer domain is an integrally closed domain such that $d=t$ \cite[Proposition 34.12]{g1}
and a P$v$MD is an integrally closed domain such that $w=t$ \cite[Theorem 3.1]{K}.

The $v$-, $t$- and $w$-operations on $R$ can be extended to the set  of nonzero
$R$-submodules of $K$ by setting
$E^v:= (R:(R:E))$,
 $E_t=\bigcup\{F_v;\, F\sub E \,; F
\mbox{ finitely generated}\}$ and $E_w=\bigcap\{E_M; \,
M\in\tmax(R)\}$, for each non zero $R$-submodule $E$ of $K$. In this way, one obtains
 \emph{semistar operations} on $R$.
 For more details, see for example \cite{FL}.
     By viewing $w$ as a semistar operation on $R$, we can say that an overring $T$ of a domain
$R$ is \emph{$t$-linked} over $R$ if $T_w=T$ \cite[Proposition 2.13]{DHLZ}. Each overring of $R$ is $t$-linked precisely when $d=w$ \cite[Theorem 2.6]{DHLZ}.

We denote by $\tspec(R)$
the set of $t$-prime ideals of $R$. Each height-one prime is a $t$-prime
and each prime minimal over a $t$-ideal is a $t$-prime.
We say that $R$ has $t$-{\it dimension
one} if each $t$-prime ideal has height one.

\medskip
We now define the ring-theoretic properties considered in this paper.

\smallskip
 \emph{The $t\#\#$-property.}
The $\#$-property and the $\#\#$-property were  introduced by R. Gilmer  \cite{g3} and R. Gilmer and W. Heinzer \cite{gh} respectively. A  domain  $R$  has the \emph{$\#$-property}, or it is a   \emph{$\#$-domain}, if $\bigcap_{M \in \mc M_1} R_M \ne
\bigcap_{M \in \mc M_2} R_M$, for any pair of distinct nonempty sets $\mc
M_1$ and $\mc M_2$ of  maximal ideals. Any overring of a one-dimensional Pr\"ufer $\#$-domain is a $\#$-domain \cite[Corollary 2]{g3}, but in general the $\#$-property is not inherited by  overrings \cite[Section 2]{gh}.  One says that $R$  has the \emph{$\#\#$-property}, or it is a   \emph{$\#\#$-domain}, if each overring of  $R$  is a   $\#$-domain.

The $t\#$-property was introduced and studied in \cite{GHL}.
A domain $R$ has the \emph{$t\#$-property} (or is a
\emph{$t\#$-domain}) if $\bigcap_{M \in \mc M_1} R_M \ne
\bigcap_{M \in \mc M_2} R_M$ for any two distinct subsets $\mc
M_1$ and $\mc M_2$ of $\tmax(R)$. Although in  \cite{GHL} it was explored
the transfer of the $t\#$-property to some distinguished classes of overrings,
it was not given any definition for the $t\#\#$-property.
Here, we say that $R$ has
the \emph{$t\#\#$-property} (or is a \emph{$t\#\#$-domain}) if
$\bigcap_{P \in \mc P_1} R_P \ne \bigcap_{P \in \mc P_2} R_P$ for
any two distinct subsets $\mc P_1$ and $\mc P_2$ of pairwise
incomparable $t$-prime ideals of $R$.
Our definition is motivated by the fact that for a P$v$MD
this is equivalent to say that
each $t$-linked overring has the $t\#$-property
\cite[Proposition 2.10]{GHL}.

\smallskip
\emph{The $tRTP$-property.} If  $R$ is an integral domain and $M$  is a unitary $R$-module, the \emph{trace} of  $M$  is the ideal of $R$ generated by the set
$\left\{f(m)\,; f \in \op{Hom}_R(M, R)\,, m \in M\right\}$.
An ideal  $J$ of $R$ is called a \emph{trace ideal}  or a \emph{strong ideal} if it is the trace of some  $R$-module  $M$.
This happens if and only if  $J = I(R:I)$, for some nonzero ideal $I$ of $R$, equivalently  $(J:J) = (R:J)$ \cite[Lemmas 4.2.2. and 4.2.3]{fhp}.
If  $V$  is a valuation domain, a trace ideal is either equal to  $V$  or it is prime \cite[Proposition 4.2.1]{fhp};  this fact led  to the consideration of several conditions related to trace ideals \cite{hp2}. The radical trace property was introduced by W. Heinzer and I. Papick \cite{hp} and is particularly significant for Pr\"ufer domains \cite{hp, l}.
$R$ is a domain satisfying the  \emph{radical trace property}, or it is an \emph{$RTP$-domain}, if  each proper strong ideal is a radical ideal, that is, for each nonzero ideal $I$ of $R$, either $I(R:I)=R$ or $I(R:I)$ is a radical ideal.

A. Mimouni studied trace properties in the setting of star operations, in particular he considered the $t$-operation \cite{M}. As in \cite{M}, we say that
a domain $R$ has the \emph{$t$-radical trace property},
or it is a {\it $tRTP$-domain}, if each proper strong
$t$-ideal of $R$ is a radical ideal. This is equivalent to say that,  for each nonzero ideal $I$ of $R$, either $(I(R:I))_t=R$ or $(I(R:I))_t$ is a radical ideal.

\smallskip
 \emph{$w$-divisoriality.}  The class of domains in which each
 nonzero ideal is divisorial has been studied, independently and with
 different methods, by H. Bass
 \cite{B}, E. Matlis \cite{Ma1} and W. Heinzer \cite{H} in the sixties.
 Following S. Bazzoni
and L. Salce \cite{BS}, a domain   in which each
 nonzero ideal is divisorial is now called a {\it divisorial domain}  and a domain such that each overring is divisorial is called  \emph{totally divisorial}.

The most suitable  star analogue of divisoriality is the notion of $w$-divisoriality that
was introduced and
extensively studied in \cite{ss}.
A \emph{$w$-divisorial domain} is a domain such that each $w$-ideal is
divisorial.

\smallskip
\emph{$w$-stability.}
Motivated by earlier work of H. Bass \cite{B} an J. Lipman
\cite{Lip} on the number of generators of an ideal, in 1974
 J. Sally and W. Vasconcelos defined a Noetherian ring $R$ to be
\emph{stable} if each nonzero ideal of $R$ is projective over its
endomorphism ring End$_R(I)$ \cite{SV1, SV2}.  When $I$ is a nonzero ideal of a domain $R$,
then End$_R(I)=(I:I)$; thus a domain $R$ is stable if each nonzero
ideal $I$ of $R$ is invertible in the overring $(I:I)$.
 B. Olberding showed that stability and divisoriality are strictly connected and
that stability is particularly significant in the context of
Pr\"ufer domains \cite{fhp, O, O3, O4,  O2}.

We introduce the notion of $w$-stability in Section 2.  We say that a $w$-ideal  $I$ of a domain $R$ is \emph{$w$-stable} if $I$
is $w$-invertible in the overring $E(I):=(I:I)$, that is $(I(E(I):I))_w=E(I)$, and say that $R$ is  \emph{$w$-stable} if each $w$-ideal  of $R$ is $w$-stable.
For a more general notion of stability with respect to a semistar operation we refer the reader to the forthcoming paper \cite{GP}.

 \section{$t\#\#$-Property and $t$-radical trace property }

    The $\#\#$-property and the radical trace property are closely
related for a Pr\"ufer domain. In this section we compare the
$t$-analogues of these two properties for P$v$MDs.

Several characterizations of P$v$MDs satisfying the $t\#\#$-property have been given in \cite[Section 2]{GHL}. For the study of the $tRTP$-property,
we need some results on branched $t$-primes.
Recall that a prime ideal $P$ of a domain $R$ is \emph{branched} if there exists
a $P$-primary ideal distinct from $P$. Clearly $P$ is branched if and
only if $PR_{P}$ is branched.
Since the localization of a P$v$MD at a $t$-prime is a valuation domain, the branched
$t$-primes of P$v$MDs can be characterized by properties similar to
those
well known for the branched primes
of Pr\"ufer domains \cite[Theorem 23.3 (e)]{g1}.

    \begin{lemma} \label{tmin} Let $R$ be a P$v$MD and $J:=x_1R+\dots+x_nR$ a nonzero finitely generated ideal such that $J_v\neq R$. If $P$ is a $t$-prime containing $J$, then $P$ is minimal over $J_v$ if and only if $P$ is minimal over $J$,  if and only if $P$ is minimal over $x_iR$, for some $i$, $1\leq i\leq n$.
\end{lemma}
    \begin{proof} It is enough to observe that, since $R_P$ is a valuation domain, we have
$J_vR_P=J_tR_P=JR_P=x_iR_P$, for some $i$, $1\leq i\leq n$.
    \end{proof}

\begin{proposition}\label{branch}
    Let $R$ be a P$v$MD and $P$ a $t$-prime of $R$. The following
conditions are
equivalent:

    \begin{itemize}
\item[(i)] $P$ is branched;
\item[(ii)] $P$ is a minimal prime of a principal ideal;
\item[(iii)] $P$ is a minimal prime of a finitely generated ideal;
\item[(iv)] $P$ is a minimal prime of a $v$-finite divisorial ideal;
\item[(v)] $P$ is not the union of the set of ($t$)-primes of $R$
properly
contained in $P$.
\end{itemize}

\end{proposition}

\begin{proof} The equivalence of conditions $(i)$, $(ii)$ and $(v)$ is
obtained by localizing at $P$ and using \cite[Theorem 17.3
(e)]{g1}. The equivalence of conditions $(ii)$, $(iii)$ and $(iv)$
follows from Lemma \ref{tmin}.
\end{proof}

 In any commutative ring with unity, if each minimal prime of an ideal $I$ is the radical of a finitely
generated ideal, then $I$  has only
finitely many minimal primes  \cite[Theorem 1.6]{gh2}. By passing through the $t$-Nagata ring, we now show that a similar result holds for $t$-ideals of P$v$MDs.

If $R$ is an integral domain, we set $R(X):= R[X]_{N}$, where $N=\{f\in R[X]:\,
c(f)=R\}$ and  $\Na:=R[X]_{N_t}$, where $N_t=\{f\in R[X];\,
c(f)_t=R\}$. $R(X)$ is called the
\emph{Nagata ring} of $R$ and
$\Na$ the \emph{$t$-Nagata ring of $R$} \cite{FL, Kang, K}.

 B. G.   Kang proved that $R$ is a P$v$MD if and only if
$\Na$ is a Pr\"ufer (indeed a Bezout) domain \cite[Theorem 3.7]{K}. In addition,  there is a lattice isomorphism
between the lattice of
$t$-ideals of $R$ and the lattice of ideals of
$\Na$ \cite[Theorem 3.4]{K}. More precisely, we have:

\begin{proposition} \label{corr} Let $R$ be a P$v$MD. Then the map
$I_t\mapsto I\Na$ is an order-preserving
bijection between the set of $t$-ideals of $R$ and the set of nonzero
ideals of $R\langle X \rangle$, whose inverse is the map $J\mapsto
J\cap R$. Moreover, $P$ is a $t$-prime (respectively, $t$-maximal)
ideal of $R$ if and only if $PR\langle X \rangle$ is a
prime (respectively, maximal)
ideal of $R\langle X \rangle$ and we have $\Na_{P\Na}= R[X]_{PR[X]}=
R_{P}(X)$.
\end{proposition}

\begin{proposition} \label{tmin2}  Let $R$ be a P$v$MD and $I$ a proper $t$-ideal of $R$. If each minimal
prime of $I$ is the radical of a $v$-finite divisorial ideal, then $I$ has finitely many minimal $t$-primes.
\end{proposition}

\begin{proof} Each minimal prime of a $t$-ideal $I$ is a $t$-ideal of $R$. By Proposition \ref{corr},
the map
$P\mapsto PR\langle X \rangle$ is a bijection between
the set of minimal primes of $I$ and the set of minimal
primes of $IR\langle X \rangle$. Moreover if $J$ is a
nonzero finitely generated ideal of $R$ such that $P=\rad(J_v)$,
then
$PR\langle X \rangle=\rad(JR\langle X \rangle)$.
Hence each minimal prime of $I\Na$ is the radical of a finitely
generated ideal. By   \cite[Theorem 1.6]{gh2},  $I\Na$ has
finitely many minimal primes and the same holds for $I$.
\end{proof}

If $T$ is an overring of $R$, the $w$-operation and the $t$-operation on $R$, viewed as semistar operations, induce two semistar operations of finite type on $T$, which here are still denoted  by $w$ and $t$ respectively. If in addition $T$ is $t$-linked over $R$, the $w$-operation is a star operation on $T$ \cite[Proposition 3.16]{EF}.  Note that this star operation, being spectral and
of finite type  \cite{FL}, is generally smaller than the
$w$-operation on $T$, that we denote by $w'$ to avoid
confusion.

\begin{proposition}\label{tlink} Let $R$ be a
P$v$MD and $T$  a $t$-linked overring of $R$. Then $T$ is a P$v$MD and $w=t=t'=w'$ on $T$, where $w'$ and $t'$ denote respectively the $w$-operation and the $t$-operation on $T$.
\end{proposition}

\begin{proof}  When $R$ is a P$v$MD also $T$ is
a P$v$MD \cite[Theorem 3.8 and Corollary 3.9]{K}. In addition, if
$R$ is a P$v$MD the two semistar operations $w$ and $t$ coincide
\cite[Theorem 3.1 ($(i)\ra(vi)$)]{FJS}. Hence  $w=t$ and $w'=t'$
as star operations on $T$. 

We next show that $t=t'$ on $T$. Let
$I$ be a nonzero ideal of $T$. Clearly, $I_t\sub I_{t'}$. On the
other hand,  we have $I_{t'} = \bigcap \{IT_M \,;\, M \in
t'\op{-Max}(T)\}$. Since $I_t= \bigcap \{I_t T_N \,;\, N \in
t\op{-Max}(T)\}$ \cite[Proposition 4]{gr}, to show that
$I_{t'}\sub I_{t}$ it suffices to show that $t\op{-Max}(T)\sub
t'\op{-Max}(T)$.

If $N\in t\op{-Max}(T)$, we have that $(N\cap R)_t \sub N_t\cap
R_t=N\cap R$. Hence $N\cap R$ is a $t$-prime of $R$ and $T_N
\supseteq R_{N\cap R}$ are valuation domains. It follows that $N$
is a $t'$-prime of $T$. In addition $N$ is $t'$-maximal because it
is $t$-maximal and each $t'$-prime of $T$ is also a $t$-prime.
\end{proof}

If $I$ is a $w$-ideal of $R$, then it is easily shown that $E(I)_w=E(I)$. Thus $E(I)$ is a $t$-linked overring of $R$.  It follows that, when $R$ is a P$v$MD, by Proposition \ref{tlink},
 $E(I)$ is  a P$v$MD and $w=t=t'=w'$ on $E(I)$.

\begin{proposition}\label{dual}
  Let $R$ be a P$v$MD. Then:

\begin{itemize}
    \item[(1)] If $I$ is a strong $t$-ideal of $R$, then $E(I)=S \cap T$,
where
$S=\bigcap_{P\in
\op{Min}(I)} R_{P}$ and $T:= \bigcap_{M \in
\tmax(R), M \nsupseteq I} R_{M}$.

\item[(2)] If $P$ is a $t$-prime of
$R$ which is not
$t$-invertible, then  $E(P) = (R:P) = R_P \cap T$, where $T:=
\bigcap_{M
\in\tmax(R), M \nsupseteq P} R_{M}$.

\item[(3)]  If
$P$ is a $t$-prime of
$R$ which is not $t$-invertible and $R$ satisfies the ascending
chain condition on radical $t$-ideals, then
$P$ is $t$-maximal in $E(P)$.
\end{itemize}
\end{proposition}

\begin{proof} (1) follows from \cite[Theorem 4.5]{HKLM}.

(2) If $P$ is not
$t$-invertible, then $P$ is strong by \cite[Proposition 2.3 and Lemma
1.2]{HZ}. Whence (2) follows from (1).

(3) By part (2), $E(P)=(R:P)$. Since $E(P)$ is $t$-linked over $R$, then $E(P)$ is $t$-flat on
$R$ (that is $E(P)_Q= R_{Q\cap R}$ for each $t$-prime ideal $Q$ of $E(P)$) \cite[Proposition 2.10]{KP} and $P$ is a $t$-ideal of $E(P)$ (Proposition \ref{tlink}).  Let
$Q$ be a $t$-prime of $E(P)$ properly containing $P$. By
$t$-flatness we can write $Q=(P'E(P))_t$, where $P'=Q\cap R$ is a
$t$-prime of $R$ properly containing $P$ \cite[Proposition 2.4]{elb1}. By
the ascending chain condition on radical $t$-ideals,
 $P'= \rad(J_t)$ for some finitely generated  ideal $J$ \cite[Lemma
3.7]{elb1}. Since $P\subsetneq P'$,   by checking
$t$-locally, we get that $P\subsetneq J_t$. We have $R=(J(R:J))_t
\sub (J(R:P))_t = (JE(P))_t\sub (P'E(P))_t= Q$. A contradiction.
Hence $P=Q$ and so $P$ is a $t$-maximal ideal of $E(P)$.
\end{proof}

\begin{lemma}\label{branchsharp} Let $R$ be a  P$v$MD satisfying the $tRTP$-property. If $P$
is a branched $t$-prime of
$R$ which is not $t$-invertible, then
$R_P \nsupseteq T$, where $T = \bigcap_{M \in
\tmax(R), M \nsupseteq P} R_{M}$.
\end{lemma}
 \begin{proof}
If $R_P \supseteq T$ then, by Proposition \ref{dual} (2), $T = E(P)$. Let $Q$
be a $P$-primary ideal of $R$. Since $P$ is a $t$-ideal, we may assume that $Q$ is a $t$-ideal. We have $QT \sub PT  =P \sub R$ and so
$QT \sub QR_P\cap R=Q$. Hence $QT=Q$. If $M$ is a $t$-maximal ideal
of $R$ such that $P\nsub M$, then $Q\nsub M$. Thus $(R:Q) \sub R_M$ and it
follows that $(R:Q) \sub T$. Hence $Q(R:Q)=Q$. Since  $R$ is a $tRTP$-domain, then
we must have $Q = P$. It follows that $P$ is not branched.
\end{proof}

\begin{theorem}\label{trtp} Let $R$ be a P$v$MD. The following
conditions are equivalent:

\begin{itemize}

    \item[(i)] $R$ is a $tRTP$-domain;

\item[(ii)] Each branched $t$-prime $P$ contains a finitely generated ideal $J$ such that $J \sub P$ and $J \nsubseteq M$, for each $M \in \tmax(R)$ not containing $P$;

\item[(iii)] Each branched $t$-prime is the radical of a $v$-finite
divisorial ideal;

\item[(iv)] Each nonzero principal  ideal has at most finitely many
minimal ($t$)-primes;

\item[(v)] Each nonzero finitely generated  ideal has at most finitely
many minimal $t$-primes;

\item[(vi)] Each  $v$-finite divisorial ideal has at most finitely many
minimal ($t$-)primes.

\end{itemize}

\end{theorem}

\begin{proof}
$(i) \ra (ii)$. Let $P$ be a branched $t$-prime of $R$. If $P$ is
$t$-invertible, then $P$ is $v$-finite and there is nothing to prove. If $P$ is not
$t$-invertible, then $E(P) = (R:P) = R_P \cap T$,
where $T = \bigcap_{M \in \tmax(R), M \nsupseteq P} R_{M}$ (Proposition \ref{dual} (2)). Since
$R_P \nsupseteq T$ (Lemma \ref{branchsharp}), there
exists a finitely generated ideal $J$ such that $J \sub P$ and $J
\nsubseteq
M$, for each $M \in
\tmax(R)$ not containing $P$ \cite[Lemma 3.6]{elb2}.

$(ii)\ra (iii)$ Let $P$ be a branched $t$-prime of $R$ and $J$ as in the hypothesis.
By Proposition
\ref{branch}, $P$ is minimal over a finitely generated
ideal $H$. Hence $P$ is the radical of the $v$-finite divisorial
ideal $(J+H)_v$.

$(iii) \ra (iv)$. Let $x\in R$ be a nonzero nonunit and let
$\{P_\al\}$
be the set of  minimal primes of $xR$. By Proposition \ref{branch} each
${P_\al}$ is branched. Hence by hypothesis each ${P_\al}$ is the radical of a
$v$-finite divisorial ideal. It follows from Proposition \ref{tmin2}
that $\{P_\al\}$ is a finite set.

$(iv) \lra (v)\lra (vi)$ follow from Lemma \ref{tmin}.

$(iv)\ra (iii)$. Let $P$ be a branched $t$-prime.
By Proposition \ref{branch}, $P$ is minimal over a
$v$-finite divisorial ideal $I$. Since $I$ has finitely many minimal
primes, then $P$ is the radical of a $v$-finite divisorial
ideal  by
\cite[Lemma 2.13]{GHL}.

$(iii)\ra(i)$. By \cite[Theorem 15]{M} it is enough to show that for
each
strong $t$-ideal $I$ and each minimal prime $P$ of $I$ we have
$IR_{P}= PR_{P}$.

Assume that $IR_{P} \subsetneq PR_{P}$. Then $P$ is branched,
because $R_{P}$ is a valuation domain. Thus
$P=\rad(H_v)$ for some  finitely generated ideal $H$. Let $a \in P$
be
such that $IR_P \subsetneq aR_P \sub PR_{P}$ and set $J = H+aR$. Then
$P$
is the radical of $J_v$. By checking
$t$-locally, we have that
$I \sub J_v$. In fact, let $M \in \tmax(R)$. If $P \nsubseteq M$, then
$IR_M \sub R_M =J_vR_M$. If $P \sub M$, then $IR_M \sub IR_P$. Hence
$a
\notin IR_M$ and $IR_M \subsetneq J_vR_M$. Since $I$ is strong, by
Proposition
\ref{dual}(1), $(R:I)=E(I) \sub R_P$. Whence $J(R:J)\sub P(R:I) \sub
PR_P$ and so  $J(R:J) \sub P$. A contradiction because $J$ is
$t$-invertible.
\end{proof}

Since in a Pr\"ufer domain the $t$-operation is trivial, we get the
following corollary, due to T. Lucas. The equivalence $(i)\lra(ii)$ is
\cite[Theorem 23]{l},  while $(i)\lra (iv)$ is, to our knowledge, unpublished.

\begin{corollary}\label{pruferrtp} Let $R$ be a Pr\"ufer domain. The following conditions are equivalent:
\begin{itemize}
    \item[(i)] $R$ is a $RTP$-domain;
\item[(ii)] Each branched prime is the radical of a finitely
generated ideal;

\item[(iii)] Each principal ideal has  at most finitely many
minimal primes;

\item[(iv)] Each  finitely generated  ideal has at most finitely many
minimal primes.
\end{itemize}
\end{corollary}

The following theorem was stated for Pr\"ufer domains in
\cite[Corollaries 25 and 26]{l}.

\begin{theorem}\label{sharp/rtp} Let $R$ be a P$v$MD.
    \begin{itemize}
\item[(1)] If $R$ is a $t\#\#$-domain, then $R$ is a $tRTP$-domain.

\item[(2)] If $R$ is a $tRTP$-domain and each $t$-prime is branched,
then $R$ is a $t\#\#$-domain.
\end{itemize}
\end{theorem}

\begin{proof} A P$v$MD $R$ has the $t\#\#$-property if and only if, for each $t$-prime ideal $P$,
there exists a finitely generated ideal $J\sub P$ such that
each $t$-maximal ideal containing $J$ must contain $P$
\cite[Proposition 2.8]{GHL}. Hence we can apply
Theorem \ref{trtp}, $(i)\lra(ii)$.
\end{proof}

\begin{theorem}\label{accr} Let $R$ be a P$v$MD. The following
conditions are  equivalent:

\begin{itemize}
    \item[(i)] $R$ satisfies the ascending chain condition on radical
    $t$-ideals;

    \item[(ii)] $R$ is a $tRTP$-domain satisfying the  ascending chain
condition on
prime $t$-ideals;

\item[(iii)] $R$ is a $t\#\#$-domain satisfying the  ascending chain
condition on prime $t$-ideals;

\item[(iv)] $R$ is a $t\#\#$-domain and each $t$-prime is branched.
\end{itemize}
\end{theorem}

\begin{proof} $(i) \lra (iii) \lra (iv)$  by \cite[Proposition 2.14]{GHL}.

$(ii) \lra (iii)$ By Proposition \ref {branch}, the  ascending chain
condition on prime $t$-ideals implies that each $t$-prime of $R$
is branched. Hence we can apply Theorem \ref{sharp/rtp}.
\end{proof}

Recall that a domain $R$ has \emph{finite character} (respectively, \emph{$t$-finite character}) if each nonzero element of $R$ belongs to
at most finitely many maximal (respectively, $t$-maximal) ideals.  A domain with finite character such that each nonzero prime ideal is contained
in a unique maximal ideal was called by E. Matlis  an {\it $h$-local domain}.
 Following  \cite{AZ}, we say that $R$ is a \emph{weakly Matlis
domain} if $R$ has
$t$-finite character and each $t$-prime ideal is contained
in a unique $t$-maximal ideal.

\begin{theorem}\label{WM} Let $R$ be a P$v$MD and consider the
following
conditions:
\begin{itemize}

    \item[(i)]  $R$ is a weakly Matlis domain;

    \item[(ii)] $R$ has $t$-finite character;

    \item[(iii)] $R$ has the $t\#\#$-property;

    \item [(iv)] $R$ is a $tRTP$-domain.

\end{itemize}

Then    $(i) \ra (ii) \ra (iii)\ra(iv)$.

If in addition each $t$-prime ideal of $R$ is
contained in a unique $t$-maximal ideal, all these conditions are
    equivalent.
    \end{theorem}

\begin{proof}  $(i)\ra(ii)$ by definition.

    $(ii) \ra (iii)$. If $R$ has $t$-finite character, for
each $\La \sub \tmax(R)$ the multiplicative
system of ideals  $\mc F(\La):=\{I \,;\, I \nsubseteq M, \text{for
each} \,M\in
\La\}$ is finitely generated
\cite[Proposition 2.7]{gab2}.  We conclude that  $R$ is a
$t\#\#$-domain
by applying
\cite[Proposition  2.8]{GHL}.

$(iii) \ra (iv)$ By Theorem \ref{sharp/rtp} (1).

\smallskip
Now assume that each $t$-prime of $R$ is contained in a unique
$t$-maximal ideal. Then clearly conditions $(i)$ and $(ii)$ are
equivalent.

$(iv)\ra(ii)$ By Theorem \ref{trtp}, for each nonzero nonunit $x\in
R$, the ideal $xR$ has finitely many minimal ($t$)-primes. Since
each $t$-prime is contained in a unique $t$-maximal ideal, then
$x$ is contained in finitely many $t$-maximal ideals.
\end{proof}

When $R$ is a Pr\"ufer domain, for $d=t$, from Theorem \ref{accr} we get \cite[Theorem 2.7]{hp} and
from Theorem \ref{WM} we get
\cite[Proposition 3.4]{O}.

\begin{remark} \rm The hypothesis that $R$ is a P$v$MD in Theorems
 \ref{accr} and \ref{WM} cannot be relaxed. In fact each Noetherian
domain is a $t\#\#$-domain \cite[Proposition 2.4]{GHL}, but it is
not necessarily a $tRTP$-domain \cite[Corollary 2.2]{hp}.
\end{remark}

A \emph{strongly discrete valuation domain} is a valuation domain such that each nonzero prime ideal  is not idempotent \cite[p. 145]{fhp} and a  \emph{strongly discrete Pr\"{u}fer domain} is a domain whose localizations at nonzero prime ideals are  strongly discrete valuation domains; equivalently  a domain such that $P\neq P^2$ for each nonzero prime ideal $P$ \cite[Proposition 5.3.5]{fhp}.
We say that a $PvMD$  $R$ is {\it strongly discrete} if $R_P$ is a strongly discrete valuation
domain for each $t$-prime ideal $P$ of $R$; equivalently, if  $(P^2)_t\not=P$, for each $P\in\tspec(R)$  \cite[Lemma 3.4]{ss}.
\emph{Generalized Krull domains} were introduced by the first author in \cite{elb1} and can be defined as  strongly
discrete P$v$MDs satisfying the ascending chain condition on
radical $t$-ideals \cite[Theorem 3.5 and Lemma 3.7]{elb1}. In the
Pr\"ufer case, that is for $d=t$, this class of domains coincides
with the class of \emph{generalized Dedekind domains} introduced by N.
Popescu  in \cite{p}.  A Krull domain is a generalized
Krull domain of $t$-dimension one \cite[Theorem 3.11]{elb1}.

\begin{theorem}\label{GK1}  Let $R$ be an integral domain. The following conditions are
equivalent:

\begin{itemize}
    \item[(i)] $R$ is a generalized Krull domain;

    \item[(ii)] $R$ is a strongly discrete P$v$MD satisfying the
$t\#\#$-property;

    \item[(iii)]  $R$ is a strongly discrete P$v$MD satisfying the
$tRTP$-property.

\end{itemize}
    \end{theorem}

    \begin{proof} $(i)\ra(ii)$ follows from Theorem
\ref{accr} and $(ii)\ra(iii)$ follows from Theorem \ref{sharp/rtp}(1).

$(iii)\ra (i)$. By Theorem \ref{trtp}, $R$ is a strongly discrete P$v$MD such that each proper $v$-finite divisorial ideal has finitely many minimal primes. Hence $R$ is a generalized Krull domain by  \cite[Theorem 3.9]{elb1}.
        \end{proof}

In the Pr\"ufer case,  for $d=t$,  we recover  from Theorem \ref{GK1}  a well known characterization of  generalized Dedekind
 domains, see for example \cite[Theorem 5.5.4]{fhp}.

\begin{remark} \rm Any $w$-divisorial domain is a $t\#$-domain. In
fact,
clearly all the $t$-maximal ideals of a $w$-divisorial domain are
divisorial; hence we can apply \cite[Theorem 1.2]{GHL}.

If $R$ is a domain such that $R_{\mc
F(\La)}:=\bigcap_{P\in\La} R_P$ is
$w$-divisorial, for each set
$\La$ of pairwise incomparable
$t$-primes, then $\mc F(\La)$ is
$v$-finite by
\cite[Proposition 2.2]{ss}. Thus $\tmax(R_{\mc
 F(\La)})=\{P_{\mc F(\La)}\,;\, P\in\La\}$ \cite[Lemma
2.1]{ss}. It follows that,
given two different sets
$\La_1$ and $\La_2$  of pairwise incomparable
$t$-primes of $R$, we have  $R_{\mc F(\La_1)}\neq
R_{\mc F(\La_2)}$. Therefore $R$ is a $t\#\#$-domain.

Conversely, it is not true that a $t\#\#$-domain is $w$-divisorial. In
fact each Noetherian domain has the $t\#\#$-property \cite[Proposition
2.4]{GHL}, but a $w$-divisorial Noetherian domain must have
$t$-dimension
one \cite[Theorem 4.2]{ss}.
\end{remark}

An integral domain $R$ is an \emph{almost Krull
domain} if
$R_M$ is a rank-one discrete valuation domain for each $t$-maximal
ideal $M$ of
$R$. Almost Krull domains
were studied by Kang
under the name
of $t$-almost Dedekind domains in \cite[Section IV]{K}. A Krull domain is an almost Krull
 domain with $t$-finite character. In dimension one,
the class of almost Krull domains coincides with the class of
almost Dedekind domains introduced by R. Gilmer \cite {g2}.
Gilmer showed that an almost Dedekind domain
satisfying the $\#$-property must be Dedekind \cite[Theorem 3]{g3}. Next, we extend this
result  to almost Krull domains.
First, we give the following characterization of  almost Krull domains, which follows directly from the definitions.

\begin{lemma}\label{akrull1} Let $R$ be an integral domain. Then $R$
is an almost Krull domain if and only if $R$ is a strongly discrete
P$v$MD of $t$-dimension one.
\end{lemma}

\begin{theorem}\label{akrull2}  Let $R$ be an integral domain. Then
the following conditions are equivalent:

\begin{itemize}

    \item [(i)] $R$ is an almost Krull domain satisfying the
$t\#$-property;

\item[(ii)] $R$ is an almost Krull domain satisfying the
$t\#\#$-property;

\item[(iii)]  $R$ is a Krull domain.

    \end{itemize}

    \end{theorem}

\begin{proof}

$(i)\ra (ii)$ Since an almost Krull domain has $t$-dimension one
(Lemma \ref{akrull1}).

$(ii)\ra (iii)$ By Lemma \ref{akrull1} and Theorem \ref{GK1}, if
(ii) holds, $R$ is a generalized Krull domain of $t$-dimension 1.
Hence $R$ is a Krull domain by \cite[Theorem 3.11]{elb1}.

$(iii)\ra (i)$ Follows from Theorem \ref{GK1}.
\end{proof}

We end this Section by putting into evidence that the $t\#\#$-property and the $tRTP$-property of a P$v$MD are related respectively to the $\#\#$-property and the $RTP$-property of its $t$-Nagata ring.

\begin{theorem}\label{Nagata1} Let $R$ be a P$v$MD. Then:
\begin{itemize}
    \item[(1)] $R$ is a $t\#$-domain if and only if
$R\langle X \rangle$ is a $\#$-domain.

        \item[(2)] $R$ is a $t\#\#$-domain if and only if
$R\langle X \rangle$ is a $\#\#$-domain.

    \item[(3)]  $R$ is a $tRTP$-domain if and only if
$R\langle X \rangle$ is a RTP-domain.
\end{itemize}

\end{theorem}

\begin{proof} (1) follows from \cite[Theorem 3.6]{GHL}.
The proof of (2) is similar and it is obtained by using Proposition \ref{corr} and the characterization of
Pr\"ufer $\#\#$-domains and P$v$MDs satisfying the $t\#\#$-property proved respectively in  \cite[Theorem 3]{gh}
and  \cite[Proposition 2.8 (7)]{GHL}.

(3) Follows from Proposition \ref{corr},Theorem \ref{trtp} $(v)$ and
Corollary \ref{pruferrtp} $(iv)$.
\end{proof}

\section{$w$-Divisoriality and $w$-stability}

The notion of $w$-divisoriality has been studied in \cite{ss}.  A domain $R$ is an integrally closed $w$-divisorial domain if and only if it is a weakly Matlis P$v$MD such that each $t$-maximal ideal is $t$-invertible \cite[Theorem 3.3]{ss} and $R$ is an integrally closed domain such that each $t$-linked overring is $w$-divisorial if and only if it is a weakly Matlis strongly discrete P$v$MD, equivalently $R$ is a $w$-divisorial generalized Krull domain \cite[Theorem 3.5]{ss}.  We now introduce the notion of $w$-stability and show that in P$v$MDs $w$-divisoriality and $w$-stability are strictly related; thus extending some results proved by B. Olberding for Pr\"ufer domains.

As before, if $T$ is a $t$-linked overring of $R$, we denote by $w$ the star operation induced on $T$ by the $w$-operation on $R$ and by $w'$ the $w$-operation on $T$. 
We say that a $w$-ideal  $I$ of $R$ is \emph{$w$-stable} if $I$
is $w$-invertible in the ($t$-linked) overring $E(I):=(I:I)$, that is if $(I(E(I):I))_w=E(I)$, and we say that $R$ is  \emph{$w$-stable} if each $w$-ideal  of $R$ is $w$-stable. 

Our first result is a generalization of \cite[Theorems 3.3 and 3.5]{O2} and shows in particular that the study of $w$-divisorial domains can be reduced to the $t$-local case.  We recall that a valuation domain is stable if and only if it is strongly discrete \cite[Proposition 5.3.8]{fhp}.

\begin{lemma}\label{qlstability} Let $R$ be a quasi-local domain.
Then  a nonzero ideal $I$ of $R$ is stable if and only if $I^2=xI$
for some $x\in I$.
\end{lemma}

\begin{proof} This follows from \cite[Lemma 3.1]{O2} and \cite[Lemma 7.3.4]{fhp}.
\end{proof}

\begin{theorem}\label{wstability} Let $R$ be an integral domain. The following conditions are
equivalent:
\begin{itemize}
\item[(i)]
$R$ is $w$-stable;
\item[(ii)] Each $w$-ideal $I$ of $R$ is divisorial in $E(I)$;
\item[(iii)] $R$ has $t$-finite character and $R_M$ is stable for each $t$-maximal ideal
$M$ of $R$.
\end{itemize}
\end{theorem}

\begin{proof}
$(i)\ra(ii)$. Let $I$ be a $w$-ideal of $R$. Denote by
$v'$ the $v$-operation on $E(I)$. Since $I$ is $w$-stable,
we have $E(I)=(I(E(I):I))_w\sub(I_{v'}(E(I):I))_w\sub E(I)$. Hence $(I_{v'}(E(I):I))_w=E(I)$ and
$I=IE(I)=(I(E(I):I)I_{v'})_w=((I(E(I):I))_wI_{v'})_w=I_{v'}$.

 $(ii)\ra(i)$. Let $I$ be a $w$-ideal of $R$ and  set $J=(E(I):I)$. Proceeding like in the proof of \cite[Theorem 3.5 ($(2)\ra(1)$)]{O2}, we have
$(E(I): IJ)=E(I)$ and hence $E((IJ)_w)=E(I)$. Thus $(IJ)_w$ is a
divisorial ideal of $E(I)$.
It follows that $(IJ)_w=(E(I):(E(I):IJ))= E(I):E(I)=E(I)$, that is $I$ is a $w$-stable ideal.

$(i)\ra(iii)$ Let $M$ be a $t$-maximal ideal of $R$ and let
$I=JR_M$ be a nonzero ideal of $R_M$, where $J$ is an ideal of $R$
which can be assumed to be a $w$-ideal (since
$J_wR_M=JR_M$). By $w$-stability in $R$, $(J(E(J):J))_w=E(J)$; in
particular, $J(E(J):J)R_M=E(J)R_M$. Since $1\in
E(J)R_M=J(E(J):J)R_M\sub I(E(I):I)\sub E(I)$, then
$I(E(I):I)=E(I)$. Hence $I$ is a stable ideal of $R_M$ and
therefore $R_M$ is stable.

We next show that $R$ has $t$-finite character. Let $M$ be a
$t$-maximal ideal of $R$. Since $R_M$ is a quasi-local stable
domain, then $M^2R_M=mMR_M$ for some $m\in M$  (Lemma
\ref{qlstability}). Set $I(M):=mR_M\bigcap R$. The ideal $I(M)$ is a
$t$-ideal of $R$ and $M^2\sub I(M)$. Hence $M$ is the only
$t$-maximal ideal of $R$ containing $I(M)$. From this, and by
checking $t$-locally, we get $(I(M):I(M))=R$. Since $R$ is $w$-stable, $I(M)$ is a
$w$-invertible ideal of $R$. Thus $I(M)$ is divisorial.

Now, let $\{M_\al\}$ be a family of $t$-maximal ideals of $R$ such
that $\bigcap M_\al\not=(0)$. We want to show that $\{M_\al\}$ is a finite family.

Set
$I_\al:=I(M_\al)$
  and let $J_\al:=(\Sigma_{\be\not=\al}(R: I_\be))_w$.   Note that $J_\al$ is a
  fractional ideal of $R$ since $\bigcap_{\be\not=\al}I_\be\supseteq
  \bigcap_{\be\not=\al}M_\be^2\not=0$.
We claim that $(J_\al:J_\al)=R$. To show this, we 
   prove that $(J_\al:J_\al)\sub R_M$ for each $t$-maximal ideal $M$ of $R$.
   Let $x\in (J_\al:J_\al)$. We first assume that
   $M\notin\{M_\be\}_{\be\neq\al}$. Since $M_\be$ is the only $t$-maximal ideal of $R$
   containing $I_\be$, then $(R:I_\be)\sub R_M$ for each $\be\not=\al$. Hence
    $x\in xJ_\al\sub J_\al\sub R_M$. If $M=M_{\gamma}$ for some $\gamma\not=\al$. We have
    $x(R:I_{\gamma})\sub(\Sigma_{\be\not=\al} (R:I_\be))_w$, and since $I_{\gamma}$ is
    $w$-invertible, then $x\in (\Sigma_{\be\not=\al}(I_{\gamma}(R:I_\be))_w)_w$. Moreover, for
    $\be=\gamma$, we have $(I_{\gamma} (R:I_{\gamma}))_w=R$, and for $\be\not=\gamma$,
    $(R:I_\be)\sub R_{M_{\gamma}}$. Hence
     $(\Sigma_{\be\not=\al}(I_{\gamma}(R:I_\be))_w)_w\sub R_{M_{\gamma}}$. Thus $x\in
     R_M$, which prove the claim.

     Now, for each $\al$, set $T_\al:=\bigcap_{\be\not=\al}M_\be$. We claim that $T_\al\nsub N$ for
each $t$-maximal ideal $N\notin\{M_\be\}_{\be\neq\al}$.  By the $w$-stability,  $J_\al$ is a $w$-invertible ideal of $(J_\al:J_\al)=R$. In particular, $J_\al$ is $w$-finite.
Thus $(R:J_\al)_N=(R_N:{J_\al}_N)=R_N$ (since $I_\be\nsub N$ for
each $\be\not=\al$). On the other hand, we have
$(R:J_\al)_N=(R:\Sigma_{\be\neq\al}(R:I_\be))R_N
=(\bigcap_{\be\not=\al}I_\be)R_N$. Thus
$\bigcap_{\be\neq\al}I_\be\nsub N$. Since
$\bigcap_{\be\neq\al}I_\be\sub T_\al$, then $T_\al\nsub N$, in
particular, $T_\al\nsub M_\al$ for each $\al$.

Now we proceed as in the proof of \cite[Theorem 3.1]{HZ1}. Set
$T:=\Sigma T_\al$. By the above result $T$ is not contained in any
$t$-maximal ideal of $R$, hence $T_t=R$. Thus
$(\Sigma_{i=1}^nT_i)_t=R$ for some finite subset $\{T_1,\ldots,
T_n\}$ of $\{T_\al\}$. Let $\{M_1,\ldots, M_n\}$ be the
corresponding set of $t$-maximal
 ideals.
If $M_\al\notin\{M_1,\ldots,
M_n\}$ for some $\al$, then $\Sigma_{i=1}^nT_i\sub M_\al$, which is
impossible. Hence  $\{M_\al\}$ is finite. Therefore $R$ has
$t$-finite character.

$(iii)\ra(i)$ Let $I$ be a $w$-ideal of $R$ and let $M_1,\ldots,
M_n$  be the $t$-maximal ideals of $R$ containing $I$.
Since $I$ is $t$-locally stable then $IR_{M_i}=J_iE(I_{M_i})$  for
some finitely generated ideal $J_i\sub I$,  $i=1, \ldots,
n$. Choose $y\in I$ such that $y\notin M$ for each $t$-maximal
ideal $M\not=M_i$ containing the ideal $H:=\Sigma J_i$ and consider the ideal
$J:=H+Ry$ of $R$. Clearly $J$ is finitely generated. One can easily check
that $IR_N=JE(I_N)$ for each $t$-maximal ideal $N$ of $R$. We next
show that $E(I_N)=E(I)_N$ for each $t$-maximal ideal $N$ of $R$.
Let $x\in E(I_N)$. Since $I_N=JE(I_N)$, then $xJ\sub I_N$. Hence
$sxJ\sub I$ for some $s\in R\sm N$.  Let $M$ be a $t$-maximal
ideal of $R$. Then $sxI_M=sxJE(I_M)\sub IE(I_M)\sub I_M$. Thus
$sxI_M\sub I_M$ for each $t$-maximal ideal $M$ of $R$, so that
$sxI\sub I$. Hence $x\in E(I)_N$. It follows that $E(I_N)=E(I)_N$,
for each $t$-maximal ideal $N$,  and
$I=\bigcap_N I_N=\bigcap_N JE(I_N) =\bigcap_NJE(I)_N=(JE(I))_w$.

Finally,
$(I(E(I):I))_N=I_N(E(I):JE(I))_N=I_N(E(I)_N:JE(I)_N)=I_N(E(I_N):I_N)=E(I_N)=E(I)_N$,
for each $t$-maximal ideal $N$.
Therefore  $(I(E(I):I))_w=E(I)$ and so $I$ is a $w$-stable ideal of
$R$.
\end{proof}

\begin{proposition}\label{likepvmd} Let $R$ be a $w$-stable domain. Then:
\begin{itemize}

\item[(1)] Each $t$-maximal ideal of $R$ is divisorial.

\item[(2)] $\tspec(R)$ is treed.

\item[(3)] $R$ satisfies the ascending chain condition on prime $t$-ideals.

\end{itemize}
\end{proposition}

\begin{proof}
 (1) Let $M$ be a $t$-maximal ideal of $R$. If $M$ is not divisorial, then $M_v=R$.
  Thus $E(M)=(R:M)=R$ and $M$ is
$t$-invertible. Hence $M$ is divisorial, which is
impossible.

(2) and (3) follow from Theorem \ref{wstability} because
a quasi-local stable domains
 has these properties \cite [Theorem 4.11]{O2}.
\end{proof}

The previous proposition shows  that  $w$-stable domains have some properties in common with generalized Krull domains \cite{elb1}.
We now prove that $w$-stability of $t$-primes enforces
a P$v$MD to be a generalized Krull domain. For Pr\"ufer domains, this
follows from
\cite[Theorem 5]{gab1} or \cite[Theorem 4.7]{O}.

\begin{theorem}\label{gkwstable} Let $R$ be a P$v$MD. The following
conditions
are equivalent:

\begin{itemize}

\item[(i)] $R$ is a generalized Krull domain;

\item[(ii)] Each radical $t$-ideal of $R$ is divisorial and each divisorial ideal
 is $w$-stable;

 \item[(iii)] Each radical $t$-ideal of $R$ is $w$-stable;

\item[(iv)] Each $t$-prime ideal of $R$ is $w$-stable.

\end{itemize}

\end{theorem}

\begin{proof} We shall freely use Proposition \ref{tlink}.

$(i)\ra(ii)$.
  Since $\tspec(R)$ is treed and a $t$-ideal of $R$ has finitely many minimal primes,
 a radical $t$-ideal of $R$ is a $t$-product of finitely many $t$-primes
\cite[Lemma 2.5]{elb2}. Hence each radical $t$-ideal of $R$ is
divisorial \cite[Proposition 3.1]{elb2}.

 Let $I$ be a divisorial ideal of $R$. If $I$ is
$t$-invertible, hence $w$-invertible, then $E(I)=R$ and so $I$ is
$w$-stable.
 If $I$ is not $t$-invertible, then consider the ideal
$H:=(I(R:I))_w$. By \cite[Proposition 2.6]{elb2}, we have
$H=(P_1\cdots P_n)_w$, where $n\geq 1$ and  each $P_i$ is a strong
$t$-prime. Thus $E(P_i)=(R:P_i)\sub(R:H)=(R:I(R:I))=((R:(R:I)):I)=E(I)$.
Since  $P_i$
is $t$-maximal in $E(P_i)$ (Proposition \ref{dual}) and $E(P_i)$ is
$t$-linked over $R$ then $P_i$ is $t$-invertible in $E(P_i)$ by
\cite[Corollaries 3.2 and 3.6]{elb1}. Hence
$(P_i(E(P_i):P_i))_w=(P_i(E(P_i):P_i))_t=(P_i(E(P_i):P_i))_{t'}=E(P_i)$,
where $t'$ denotes the $t$-operation on $E(P_i)$. Thus $P_iE(I)$
is $w$-invertible in $E(I)$, for each $i$. It follows that $H$ is
$w$-invertible in $E(I)$ and so $I$ has the same property.

$(ii)\ra(iii)\ra (iv)$ are clear.

$(iv)\ra(i)$. Let $P$ be a $t$-prime of $R$. Since $P$ is
$w$-invertible in $E(P)$, then $P\neq (P^2)_{w}$ and so $P\neq
(P^2)_{t}$. Thus $R$ is a strongly  discrete P$v$MD.

To prove that $R$ is a generalized Krull domain, it is enough to
show that $R$ has the $t\#\#$-property (Theorem \ref{GK1}). Let
$T$ be a $t$-linked overring of $R$  and denote by $t'$ the
$t$-operation on $T$. Let $M$ be a $t'$-maximal ideal of $T$.
Since $T$ is a P$v$MD, then $T=E(M)$.  The ideal $P=M\bigcap R$ is a
$t$-prime of $R$ and $M=(PT)_{t'}=(PT)_w$ (cf. \cite[Proposition
2.10]{KP} and \cite[Proposition 2.4]{elb1}). Thus $R\sub E(P)\sub
E(M)=T$.  Since $P$ is $w$-stable, then $PT$ is
 $w$-invertible in $T$, and hence $M$ is  $w$-invertible in $T$. So,
$M$ is a $t'$-invertible $t'$-ideal of $T$ (since $w=t'$ in $T$). In
particular $M$ is a  divisorial ideal of $T$. We conclude that $T$ is
a $t\#$-domain by applying \cite[Theorem 1.2]{GHL}.
\end{proof}

It is known that a generalized Dedekind domain need not be stable \cite[Example 10]{gab1}. In fact an integrally closed domain is stable if and only if it is a strongly discrete Pr\"ufer domain with finite character \cite[Theorem 4.6]{O}. Hence a
generalized Dedekind domain is stable if and only if it has finite character.
We now extend these results to generalized Krull domains.

\begin{lemma} \label{SD=GK} A domain with $t$-finite character is a strongly discrete P$v$MD if and only if it is a generalized Krull domain.
\end{lemma}
\begin{proof} A strongly discrete P$v$MD is a generalized Krull domain if and only if
each nonzero nonunit has finitely many minimal $t$-primes \cite[Theorem 3.9]{elb1}. We conclude by recalling that in a P$v$MD two incomparable $t$-primes are $t$-comaximal.
\end{proof}

\begin{theorem} \label{wstableic} Let $R$ be an integral domain.
    The following conditions are equivalent:

\begin{enumerate}

    \item[(i)] $R$ is integrally closed and $w$-stable;

\item[(ii)] $R$ is a $w$-stable P$v$MD;

   \item[(iii)] $R$ is a strongly discrete P$v$MD with $t$-finite character;
 \item[(iv)] $R$ is a generalized Krull domain with $t$-finite character;

 \item[(v)] $R$ is a $w$-stable generalized Krull domain;

\item[(vi)] $R$ is a P$v$MD with $t$-finite character and each
$t$-prime ideal of $R$ is $w$-stable;
 \item[(vii)] $R$ is $w$-stable and each $t$-maximal ideal of $R$ is
$t$-invertible.
\end{enumerate}

\end{theorem}

\begin{proof}
    $(i)\ra(ii)$. If $M$ is a $t$-maximal ideal of $R$,
then $R_M$ is an integrally closed stable domain by Theorem \ref{wstability}.
 Hence $R_M$ is a
valuation domain \cite[Theorem 4.6]{O}.

    $(ii)\ra(iii)$. For each $t$-maximal ideal $M$ of $R$,
$R_M$ is a valuation stable domain (Theorem \ref{wstability}). Hence $R_M$ is a
strongly discrete valuation domain \cite[Proposition 5.3.8]{fhp}. The $t$-finite character follows again
from Theorem \ref{wstability}.

$(iii)\lra(iv)$ by Lemma \ref{SD=GK}.

$(iv)\ra(v)$  $R_M$ is stable, for each $M\in \tmax(R)$, because
it is a strongly discrete valuation domain
\cite[Proposition 5.3.8]{fhp}. By the $t$-finite character, $R$
is $w$-stable (Theorem \ref{wstability}).

 $(v)\ra(vii)$ because each $t$-maximal ideal of a generalized Krull domain  is
$t$-invertible  \cite[Corollary 3.6]{elb1}.

$(vii)\ra(i)$ By Theorem  \ref{wstability}, $R_M$ is a local
stable domain, for each $M\in \tmax(R)$. Since $M$ is
$t$-invertible, $MR_M$ is a principal ideal. Hence $R_M$ is a
valuation domain \cite[Lemma 4.5]{O2} and $R$ is integrally
closed.

$(iv)\lra(vi)$ follows from Theorem \ref{gkwstable}.
\end{proof}

 By Theorem \ref{gkwstable}, each divisorial ideal of a generalized Krull domain is $w$-stable. Hence
 a $w$-divisorial generalized Krull domain is $w$-stable. Several characterizations of
$w$-divisorial generalized Krull domains were given in
\cite[Theorem 3.5]{ss}. The following theorem says something more in terms
of $w$-stability; similar results for Pr\"ufer domains were
obtained by Olberding \cite{O, O4}.

\begin{theorem} \label{wdiv}
    Let $R$ be an integral domain. The following conditions are
equivalent:

\begin{enumerate}
    \item[(i)] $R$ is an integrally closed $w$-divisorial $w$-stable
    domain;
    \item[(ii)] $R$ is a $w$-stable $w$-divisorial P$v$MD;
    \item[(iii)] $R$ is a $w$-divisorial generalized Krull domain;
   \item[(iv)] $R$ is a weakly Matlis $w$-stable P$v$MD;
    \item[(v)] $R$ is a weakly Matlis strongly discrete P$v$MD;
 \item[(vi)] $R$ is a weakly Matlis generalized Krull domain.
\end{enumerate}
    \end{theorem}

   \begin{proof}
 $(i)\lra(ii)$ by Theorem  \ref{wstableic}.

 $(ii)\ra(iv)$ because a $w$-divisorial domain is weakly Matlis \cite[Theorem 1.5]{ss}.

$(iii)\lra(v)$ by \cite[Theorem 3.5]{ss}.

$(v)\lra(vi)$ by Lemma \ref{SD=GK}.

$(iii)+(vi)\ra(i)$ and $(iv)\lra(v)$ by Theorem \ref{wstableic}, because a weakly Matlis domain has $t$-finite character.
\end{proof}

From Theorems \ref{wstableic} and \ref{wdiv}, we immediately get:

\begin{corollary}
Let $R$ be an integrally closed $w$-stable domain. Then $R$ is $w$-divisorial if and only if each nonzero $t$-prime ideal is contained in a unique $t$-maximal ideal.
\end{corollary}

A domain is stable and divisorial if and only if it is totally divisorial \cite[Theorem 3.12]{O4}.  The following is the $t$-analogue of this result in the integrally closed case.

\begin{corollary} Let $R$ be an integral domain. The following conditions are
equivalent:

\begin{enumerate}
    \item[(i)] $R$ is a $w$-stable $w$-divisorial P$v$MD;
    \item[(ii)] $R$ is integrally closed and each $t$-linked overring of $R$ is $w$-divisorial.
\end{enumerate}
\end{corollary}
\begin{proof}
By \cite[Theorem 3.5]{ss}, if $R$ is integrally closed, each $t$-linked overring of $R$ is $w$-divisorial if and only if $R$ is a weakly Matlis strongly discrete P$v$MD. We conclude by applying Theorem \ref{wdiv}.
\end{proof}

We recall that each overring of a
domain $R$ is $t$-linked if and only if $d=w$ on $R$ \cite[Theorem 2.6]{DHLZ} and that each overring of a stable domain is stable \cite[Theorem 5.1]{O2}. We now prove  that $w$-stability is preserved by $t$-linked extension.

\begin{theorem}\label{wlink} Let $R$ be an integral domain and $T$ a $t$-linked overring of $R$. 
If $R$ is $w$-stable then $T$ is $w'$-stable, where $w'$ denotes the $w$-operation on $T$.
\end{theorem}

\begin{proof} We shall use Theorem \ref{wstability}. Since $R\sub T$
is $t$-linked, for each $t'$-maximal ideal $M$ of $T$, there is a
 $t$-maximal ideal $N$ of $R$ such that $R_N\sub
T_M$ \cite[Proposition 2.1]{DHLZ}. Hence $T_M$ is an
overring of a stable domain and is therefore stable \cite[Theorem
5.1]{O2}.

We next show that $T$ has $t$-finite character. Let   $N$ be a $t$-maximal ideal of $R$ and
let $\{M_\al\}$ be a family of $t$-maximal ideals of $T$ such that
$\bigcap_\al M_\al\neq(0)$ and $M_\al\cap R\sub N$ for each $\al$.
 Set $S:=\bigcap_\al T_{M_\al}$.
Then $S$ is a stable domain since it is an overring of the stable
domain $R_N$ \cite[Theorem 5.1]{O2}. The prime ideals $P_\al=
M_\al T_{M_\al}\cap S$ of $S$ are pairwise incomparable, since
$S_{P_\al}=T_{M_\al}$ for each $\al$. We have $(0)\neq\bigcap_\al
M_\al\sub\bigcap_\al P_\al$, and, since $S$ is treed \cite[Theorem 4.11
(ii)]{O2} and has finite character \cite[Theorem 3.3]{O2}, then
$\{P_\al\}$ must be finite. Hence $\{M_\al\}$ is also a
finite set. Since $R$ has $t$-finite character, it follows that
$T$ has $t$-finite character.
\end{proof}

 We  do not know whether the integral closure of a $w$-stable domain is $w'$-stable.  In fact the integral closure of a domain $R$ is not always $t$-linked over $R$
  \cite[Section 4]{DHLRZ} and  we cannot apply Theorem \ref{wlink}. However,
  the $w$-integral closure $R^{[w]}:= \bigcup\{(J_w:J_w)\,; J \mbox{ a finitely generated ideal of } R\}$  is integrally closed and $t$-linked over $R$ \cite[Proposition 2.2 (a)]{DHLZ}. Thus we immediately get:

 \begin{corollary}
 The $w$-integral closure of a $w$-stable domain is a $w'$-stable P$v$MD.
\end{corollary}

We end by remarking that, in the integrally closed case, $w$-divisoriality and $w$-stability
correspond to divisoriality and stability of the $t$-Nagata ring. We shall make use of Proposition \ref{corr}.

\begin{theorem}\label{Nagata} Let $R$ be an integral domain. Then:

\begin{itemize}

\item[(1)] $R$ has $t$-finite character if and only if
$R\langle X \rangle$ has finite character.

 \item[(2)]  $R$ is a Weakly Matlis P$v$MD if and only if
$R\langle X \rangle$ is an $h$-local Pr\"ufer domain.

\item[(3)] $R$ is a strongly discrete P$v$MD if and only if
$R\langle X \rangle$ is a strongly discrete Pr\"ufer domain.

\item[(4)] $R$ is a generalized Krull domain if and only if
$R\langle X \rangle$ is a generalized Dedekind domain.

\item[(5)] $R$ is an integrally closed $w$-divisorial domain if and
only
if $R\langle X \rangle$ is an integrally closed divisorial domain.

\item[(6)] $R$ is an integrally closed $w$-stable domain if and only
if $R\langle X \rangle$ is an integrally closed stable
domain.
\end{itemize}

\end{theorem}

\begin{proof} Denote by $c(f)$ the content of a polynomial $f(X)\in R[X]$.

(1) We have $\op{Max}(R\langle X \rangle)=\{MR\langle X \rangle ; M\in
\tmax(R)\}$.   Since $f(X)\in  MR[X]$ if and only if $c(f)_v \sub M$, if $R$  has $t$-finite character, then $R\langle X \rangle$ has $t$-finite character.
 The converse is clear.

(2) Follows from (1) and Proposition \ref{corr}.

 (3) For $M\in \tmax(R)$, we have
that
 $R\langle X \rangle_{M\Na}= R[X]_{MR[X]}= R_{M}(X)$ is a
strongly discrete
valuation domain if and only if $R_M$ has the same property.

 (4) Follows from (3) and Proposition \ref{corr} by recalling the definitions.

(5)  When $R$ is integrally closed,  $R$ is divisorial if and only if $R$
is an $h$-local Pr\"ufer domain  such that each maximal ideal is invertible \cite[Theorem 5.1]{H} and $R$ is $w$-divisorial if and only if $R$
is a weakly Matlis P$v$MD such that each $t$-maximal ideal is $t$-invertible \cite[Theorem 3.3]{ss}.
Hence we can conclude by applying part (2) and recalling that, for each $M\in \tmax(R)$,
 $MR\langle X \rangle$ is invertible if and only if $M$ is $t$-invertible \cite[Theorem 2.4]{K}.

(6) Follows from \cite[Theorem 4.6]{O}, Theorem \ref{wstableic} and  statements (1) and
(3).
\end{proof}


\end{document}